\documentclass[12pt]{article}

\usepackage{amsmath,amssymb,amsthm}
\usepackage{graphicx}
\usepackage{ifpdf}
\usepackage{color}
\allowdisplaybreaks
\DeclareGraphicsExtensions{.pdf,.png,.jpg,.jpeg,.mps,.eps}

\usepackage{times}
\usepackage[all]{xypic}
\usepackage[colorlinks=false,linktocpage=true]{hyperref}
\usepackage{hyperref}
\hypersetup{colorlinks=false}

\newtheorem{remark}{Remark}[section]
\newtheorem{theorem}{Theorem}[section]
\newtheorem{corollary}{Corollary}[section]
\newtheorem{lemma}{Lemma}[section]                                                                                                              

\newtheorem{definition}{Definition}[section]
\newtheorem{proposition}{Proposition}[section]

\numberwithin{equation}{section}

\def\R{\mathbb R}

\def\C{\mathbb C}

\def\Z{\mathbb Z}

\def\H{\mathbb H}
\def\D{\mathfrak D}

\def\F{\mathcal F}

\def\tr{{\rm tr}}
\newcommand{\PSL}{{\rm PSL}}
\newcommand{\SL}{{\rm SL}}

\def\Im{{\rm Im\,}}
\def\Re{{\rm Re\,}}

\begin{document}

\title{Fundamental domains in $\PSL(2,\R)$ for Fuchsian groups}

\author{{\sc Huynh Minh Hien} \\[1ex] 
	Department of Mathematics and Statistics,\\
	Quy Nhon University, 170 An Duong Vuong, Quy Nhon, Vietnam;\\
	e-mail: huynhminhhien@qnu.edu.vn
	\date{} 
}
\maketitle

\begin{abstract} In this paper, we provide a necessary  and sufficient condition for a set in $\PSL(2,\R)$
	or in $T^1\H^2$ to be a fundamental domain for a given Fuchsian group via its respective fundamental domain in the hyperbolic plane $\H^2$. 
\end{abstract}


\section{Introduction}
Fundamental domains arise naturally in the study of group actions on topological spaces.
 The concept {\em fundamental domain} is used to describe a set in a topological space under a group action
 of which the images tessellate the whose space.
The term {\em fundamental domain} is well-known in the model of the hyperbolic plane $\H^2$ for
the action of Fuchsian groups via M\"obius transformations. If there exists a point in $\H^2$ that is not a fixed point for all elements different from the unity in a Fuchsian group $\Gamma$ then there always exists a convex and connected fundamental domain for $\Gamma$ named Dirichlet domain (see \cite{einsward,katok}). Other examples of fundamental domains are Ford domains (see  \cite{series}). {\em  Poincar\'e's polygon theorem} \cite{hubbard,katok,series} provides a fundamental domain, which is  a polygon for the Fuchsian group generated 
by the side-pairing transformations. In this case, if the polygon has finite edges (and hence it is relatively compact), the Fuchsian group
is finitely generated and the space of $\Gamma$-orbits denoted by $\Gamma\backslash\H^2$ is compact. Fundamental 
domains have several applications for the study of $\Gamma\backslash\H^2$.  If the action of $\Gamma$ has no fixed points, 
the quotient space $\Gamma\backslash\H^2$ has a Riemann surface structure that is a closed Riemann surface of genus at least $2$  
and has the hyperbolic plane $\H^2$ as the universal covering. Furthermore, 
it is well-known that any compact orientable surface 
with constant negative curvature is isometric to a factor $\Gamma\backslash \H^2$. 
If the Fuchsian group has a finite-area fundamental domain then all the fundamental domains have finite area and
have the same area. This area is defined for the measure of the quotient space $\Gamma\backslash\H^2$. 
In addition, the space $\Gamma\backslash\H^2$ is compact if and only if every fundamental domain in $\H^2$ for $\Gamma$ is relatively compact (see \cite{katok}). 

There exists a bijection $\Theta: T^1\H^2\to\PSL(2,\R)$. The natural Riemannian metric on $\PSL(2,\R)$ induces
a left-invariant metric function (a metric in the usual sense). The topology induced from this metric is the same as the quotient topology induced from the one in $\SL(2,\R)$. The Sasaki metric  on the unit tangent bundle $T^1\H^2$ with respect to
the hyperbolic metric on $\H^2$ makes $\Theta$ an isometry. 
This induces an isometry from  $T^1(\Gamma\backslash\H^2)$ to
$\Gamma\backslash\PSL(2,\R)$, where $\Gamma\backslash\PSL(2,\R)$ denotes the collection of  right co-sets $\Gamma g$ of  $\Gamma$ in $\PSL(2,\R)$, which is also obtained from a left action of Fuchsian group $\Gamma$
on $\PSL(2,\R)$. Furthermore, there is an action of $\Gamma$ on the unit tangent bundle $T^1\H^2$ by derivative operators  and this arises fundamental domains for $\Gamma$ in $T^1\H^2$ also. However, up to now there have not been any results about fundamental domains for $\Gamma$ in $\PSL(2,\R)$ or in $T^1\H^2$.
The aim of this paper is to study fundamental domains in $\PSL(2,\R)$ and in $T^1\H^2$ for Fuchsian groups. A necessary and sufficient condition
for a set in $\PSL(2,\R)$ or in $T^1\H^2$ to be a fundamental domain via its respective fundamental domain in $\H^2$ is provided. 

The paper is organized as follows. In the next section we present the actions of Fuchsian groups 
on the hyperbolic plane $\H^2$, the unit tangent bundle $T^1\H^2$ and  the group $\PSL(2,\R)$. The main results are stated and proved in Section 3.

      \setcounter{equation}{0}
\section{Preliminaries}
In this section we introduce the necessary background material which can be found in \cite{hubbard,katok,series}. 
The unity of an arbitrary group is always denoted by $e$.
\subsection{Fundamental domains}
Let $X$ be a non-empty set and  let $G$ be a group. Let $\rho: G\times X \to X$ be a (left) group action,
that is, $\rho(e,x)=x$ and $\rho(g_1,\rho(g_2,x))=\rho(g_1g_2,x)$ for all $x\in X$ and $g_1,g_2\in G$. 
For a subset $A\subset X$, define $\rho(g,A)=\{\rho(g,x), x\in A \}$.
\begin{definition}\label{funddo-def}\rm
	Let $G$ be a group and let $X$ be a topological space. Suppose that $\rho:G\times X\to X$ is group action. 
	A non-empty open set $F\subset X$ is said to be a {\em fundamental domain} for $G$, if
	\begin{itemize}
		\item[(a)] $\bigcup_{g\in G} \rho(g,\overline{F})=X$ and
		\item[(b)] $\rho(g, F)\cap F=\varnothing$ for all $g\in G\setminus\{e\}$.
	\end{itemize}
	Here $e$ is the unity of $G$ and $\overline{F}$ denotes the closure of $F$ in $X$.
\end{definition}
Due to the fact that $G$ is a group,
condition (b) is equivalent to
\[\rho(g_1, F)\cap\rho(g_2,F)=\varnothing
\quad\mbox{for all}\quad g_1, g_2\in G,
\,\,g_1\neq g_2. \]
We will introduce some examples in the next subsection.
\subsection{$\H^2$ and $\PSL(2,\R)$}
The hyperbolic plane is the upper half plane 
$\H^2=\{(x, y)\in\R^2: y>0\}$, endowed with the Riemannian metric $(g_z)_{z\in\H^2}$,
where $g_z(\xi,\eta)=\frac{\xi_1\eta_1+\xi_2\eta_2}{(\Im z)^2}$ 
for $\xi=(\xi_1,\xi_2), \eta=(\eta_1,\eta_2)\in  T_z\H^2\cong \C$. 
The group of M\"obius transformations ${\rm M\ddot ob}(\H^2)=\{z\mapsto \frac{az+b}{cz+d}\,:\, a,b,c,d\in\R\,, ad-bc =1\}$ 
can be identified with the projective group $\PSL(2,\R)=\SL(2,\R)/\{\pm E_2\}$ by means of the isomorphism
\begin{equation}\label{Phi}
\Phi\Bigg(\pm\bigg(\begin{array}{cc}a&b\\c&d\end{array}\bigg)\Bigg)=z\mapsto\frac{az+b}{cz+d},
\end{equation} 
where $\SL(2,\R)$ is the group of all real $2\times 2$ matrices with unity determinant, 
and $E_2$ denotes the unit matrix. 

Let $\Gamma$ be a Fuchsian group, which is a discrete subgroup in $\PSL(2,\R)$. We consider the action 
$\rho: \Gamma\times \H^2\to \H^2,\,\rho(\gamma,z)=\Phi(\gamma)(z)\ \mbox{for}\ (\gamma,z)\in \Gamma\times\H^2$.
The action is called free if
$\Phi(\gamma)(z)=z$ for some $z\in\H^2$ then $\gamma=e$. In this case, there always exist fundamental domains for $\Gamma$ as follows.

\begin{proposition}
	Let $\Gamma\subset {\rm PSL}(2, \R)$ be a Fuchsian group
	and take $z_0\in\H^2$ such that $z_0\neq\Phi(\gamma)(z_0)$ holds
	for all $\gamma\in\Gamma\setminus\{e\}$. Then the Dirichlet region
	\[ D_{z_0}(\Gamma)=\Big\{z\in\H^2: d_{\H^2}(z, z_0)<d_{\H^2}(z, T(z_0))
	\,\,\mbox{for all}\,\,T=\Phi(\gamma), \gamma\in\Gamma\setminus\{e\}\Big\} \]
	is a fundamental domain for $\Gamma$ which contains $z_0$.
\end{proposition}
See \cite[Lemma 11.5]{einsward} for a proof. Note that such a $z_0$ does exist if the action of $\Gamma$ on $\H^2$ is free. 

For $g=\Big\{\pm\big(\scriptsize\begin{array}{cc}a&b\\c&d\end{array} \big)\Big\}\in\PSL(2,\R)$,
the trace of $g$ is defined by \[\tr(g)=|a+d|.\]
Element $g$ is called hyperbolic if $|\tr(g)|>2$, elliptic if $|\tr(g)|<2$ and parabolic if $|\tr(g)|=2$.
It is well-known that the action of $\Gamma$ on $\H^2$ is free if and only if $\Gamma$ does not contain any elliptic elements (see \cite[Remark 4.26]{series}). 
 
 For any $g\in\PSL(2,\R)$, the cyclic group $\langle g\rangle =\{g^n:n\in\Z \}\subset \PSL(2,\R)$
 is a Fuchsian group. We will consider fundamental domains for $\langle g\rangle$ with some special classes of $g$.
 For $t\in\R$, let  
 \begin{eqnarray*}
 	A_t=\bigg(\begin{array}{cc} e^{t/2} & 0 \\
 		0 & e^{-t/2}\end{array}\bigg),
 	B_t=\bigg(\begin{array}{cc} 1 & t \\
 		0 & 1\end{array}\bigg), D_t=\bigg(\begin{array}{cc} \cos\frac{t}2 & \sin\frac{t}2 \\
 		-\sin\frac{t}2 & \cos\frac{t}2\end{array}\bigg)\in\SL(2,\R)
 \end{eqnarray*} 
 and respectively
 \[ a_t=[A_t], b_t=[B_t], 
 d_t=[D_t]\in \PSL(2,\R) .\]

\begin{proposition}\label{examfun} 
	(a) For any $t>0$, the set
	\[F_t=\{z=x+iy\in\H^2: 1<y<e^t\}\]
	is a fundamental domain in $\H^2$ for
	the Fuchsian group 
	$\langle a_t\rangle $.
	
	(b) For any $t>0$, the set 
	\[E_t=\{z=x+iy\in\H^2: 0< x< t \} \]
	is a fundamental domain in $\H^2$ for the Fuchsian group $\langle b_t\rangle$.
\end{proposition}
\noindent{\bf Proof\,:} (a) Obviously $F_t$ is open and
	\[\overline {F_t}
	=\{z=x+iy\in\H^2: 1\leq y \leq e^t \}.\]
	We have 
	$\Phi(a_{jt})
	=e^{jt}\,{\rm id}$, with ${\rm id}: \H^2\to\H^2$ denoting the identity map. Here 
	\[ \Phi(a_{jt})(\overline{F_t})=e^{jt}\,{\rm id}(\{z=x+iy\in\H^2: 1\le y\le e^t\})
	=\{z\in\H^2: e^{jt}\le y\le e^{(j+1)t}\}, \] 
	so that $\bigcup_{j\in\Z}\Phi(a_{jt})(\overline{F_t})=\H^2$ 
	and $\Phi(a_{jt})({F_t})\cap\Phi(a_{kt})({F_t})=\varnothing$ 
	for $j\neq k$. 
	
	(b) It is proved analogously to (a).
	{\hfill$\Box$}\bigskip

The collection of right co-sets $\Gamma g$ of $\Gamma$ in $\PSL(2,\R)$ denoted by $\Gamma\backslash\PSL(2,\R)$ can be also obtained by $\Gamma$-orbits of the left action
\begin{equation}\label{varrho}\varrho: \Gamma \times \PSL(2,\R)\rightarrow \PSL(2,\R),\,
\varrho(\gamma,g)=\gamma g\ \mbox{for}\ \gamma\in\Gamma,g\in\PSL(2,\R).\end{equation} 
This leads to  the concept {\em fundamental domain} in $\PSL(2,\R)$.

\begin{remark}\rm If $\F\subset \PSL(2,\R)$ is a fundamental domain for $\Gamma$ and $\gamma\in\Gamma\setminus\{e\}$,
	then $\gamma \F$ is a fundamental domain disjoint from $\F$. For, it is obvious that  $\gamma \F$ is open since $\F$ open 
	and $\overline{\gamma\F}=\gamma\overline{\F}$
	in $\PSL(2,\R)$. Therefore
	\[ \bigcup_{\gamma'\in\Gamma} \gamma'\overline{\gamma\F}
	=\bigcup_{\gamma'\in\Gamma} \gamma'(\gamma\overline{\F})
	=\bigcup_{\gamma'\in\Gamma}(\gamma'\gamma)(\overline{\F})
	=\bigcup_{\eta\in\Gamma}\eta \overline{\F}=\PSL(2,\R), \]
	and for $\gamma'\in\Gamma\setminus\{e\}$,
	\[ (\gamma\F)\cap \gamma'(\gamma\F)=\gamma \F\cap( \gamma'\gamma)\F=\varnothing \]
	by Definition \ref{funddo-def} (a), due to $\gamma\neq\gamma'\gamma$. {\hfill$\diamondsuit$}
\end{remark}

\subsection{$T^1\H^2$}

The unit tangent bundle of $\H^2$ is defined by
\begin{equation}\label{T1H2-def}
T^1\H^2=\{(z, \xi): z\in\H^2, \xi\in T_z\H^2,
{\|\xi\|}_z=g_z(\xi, \xi)^{1/2}=1\}.
\end{equation}
For $g\in {\rm PSL}(2, \R)$ we consider the derivative operator
\begin{equation*}\label{cald-def}
{\D}g: T^1\H^2\to T^1\H^2
\end{equation*}
defined as
\begin{equation*}\label{bursetal}
{\D}g(z, \xi)=(T(z), T'(z)\xi),
\end{equation*}
where $T=\Phi(g)$; recall $\Phi$ in \eqref{Phi}. Then ${\D}$ is well-defined. Explicitly, 
if $g=\scriptsize\Big\{\pm\big(\begin{array}{cc} a & b \\ c & d\end{array}\big)\Big\}$, 
then $T(z)=\frac{az+b}{cz+d}$ and $ad-bc=1$, whence
\begin{equation}\label{calDexpl}
{\D}g(z, \xi)=\Big(\frac{az+b}{cz+d},\,\frac{\xi}{(cz+d)^2}\Big).
\end{equation}

Let $\Gamma\subset\PSL(2,\R)$ be a subgroup. Consider the group action 
$$\kappa: \Gamma\times T^1\H^2\to T^1\H^2,\ \kappa(\gamma, (z,\xi))={\D }\gamma (z,\xi)
\
 \mbox{for}\ \gamma\in\Gamma, (z,\xi)\in T^1\H^2.$$
If $\Gamma=\PSL(2,\R)$ then the action is simply transitive (see \cite[Lemma 9.2]{einsward}), that is,
for given $(z,\xi), (w,\eta)\in T^1\H^2$, there exists a unique $g\in\PSL(2,\R)$ such that 
$\kappa(g,(z,\xi)) =\D g(z,\xi)=(w,\eta)$. 
In particular, we have the following property.

\begin{lemma}\label{unig}
	For each $(z,\xi)\in T^1\H^2$, there is a unique $g\in\PSL(2,\R)$ such that $\D g (i,i)=(z,\xi).$
\end{lemma}

Explicitly, if $(z,\xi)\in T^1\H^2$ then $g=\Big\{\pm\scriptsize\Big(\begin{array}{cc}a&b\\c&d \end{array}\Big)\Big\}$ 
is defined by
\begin{equation}\label{dg}\frac{ac+bd}{c^2+d^2}=\Re z, \frac{1}{c^2+d^2}=\Im z, \frac{2cd}{(c^2+d^2)^2}={\rm Re\,}\xi,
\quad \frac{d^2-c^2}{(c^2+d^2)^2}=\Im\xi.
\end{equation}
We will use these relations afterwards.

\section{Main results}
This section deals with the relation of fundamental domains for a Fuchsian group in $\H^2$, $\PSL(2,\R)$ and in $T^1\H^2$.

The main result of this paper is the following: 
\begin{theorem}\label{fdPSL}
	Let $\Gamma\subset \PSL(2,\R)$ be a Fuchsian group. For $F\subset \H^2$, denote 
	\[\F=\{g=b_x a_{\ln y}d_\theta\, :\, x+iy\in F, \theta\in [0,2\pi)  \} \subset \PSL(2,\R).\]
	Then $F$ is a fundamental domain for $\Gamma$ in $\H^2$ if and only if $\F$  
	is a fundamental domain for $\Gamma$ in $\PSL(2,\R)$. 
\end{theorem} 
\begin{remark}\rm
	Recall that if $\Gamma$ contains no elliptic elements  then there always exist fundamental domains in $\H^2$ for
	$\Gamma$ and hence fundamental domains in $\PSL(2,\R)$ do always exist. The collection of $\Gamma$-orbits of
	the action $\varrho$ (see \eqref{varrho}) denoted by $\Gamma\backslash\PSL(2,\R)=\{\Gamma g, g\in\PSL(2,\R) \}$
	is compact if and only if the quotient space $\Gamma\backslash\H^2$ is compact
	if and only if there is a relatively compact fundamental domain (in $\H^2$ or in $\PSL(2,\R)$) for $\Gamma$.
	In this case all fundamental domains of $\Gamma$ are relatively compacts. 
	For proofs of the case in $\H^2$, see \cite[Chapter 3]{katok}. {\hfill$\diamondsuit$}
\end{remark}

In order to prove Theorem \ref{fdPSL}, we need the following  factorization, which is called NAK decomposition (so-called Iwasawa decomposition).
\begin{lemma}[\cite{sieber}]\label{siber1}
	If $G= \scriptsize\Big(\begin{array}{cc} a&b\\c&d\end{array}\Big)\in\SL(2,\R)$ then
	$G=B_xA_{\ln y}d_\theta$ with
	\begin{equation}
	x=\frac{ac+bd}{c^2+d^2},\ \  y=\frac{1}{c^2+d^2},\ \  \theta=-2\arg(d+ic).
	\end{equation}
	
\end{lemma}

\begin{lemma}\label{lemnak}
	(a) If $g=[G]\in\PSL(2,\R)$ for $G= \scriptsize\Big(\begin{array}{cc} a&b\\c&d\end{array}\Big)\in\SL(2,\R)$ then 
	$g= b_x a_{\ln y}d_\theta$ with
	\begin{equation}\label{nak1}
	x=\frac{ac+bd}{c^2+d^2},\ \  y=\frac{1}{c^2+d^2},\ \  \theta=-2\arg(d+ic).
	\end{equation}
	(b) \begin{eqnarray*}
	\PSL(2,\R)&=&\{b_xa_{\ln y}d_\theta : x+iy\in\H^2,  \theta\in [0,2\pi) \}\\
	&=&
\{b_xa_{\ln y}d_\theta : x+iy\in\H^2, \theta\in\R \}.
\end{eqnarray*} 
\end{lemma}
\noindent
{\bf Proof\,:}
(a) This follows directly from Lemma \ref{siber1}. (b) According to (a), every element $g=\Big\{\pm \Big({\scriptsize\begin{array}{cc} a&b\\c&d\end{array}}\Big)\Big\}\in\PSL(2,\R)$
has the decomposition $g=b_xa_{\ln y}d_\theta$ for $x+iy\in\H^2$ and $\theta=-2\arg(d+ic)\in (-2\pi,0]$. It remains to 
verify that we can find some $\theta'\in [0,2\pi)$ such that $d_{\theta'}=d_\theta$
and as a consequence, $g=b_xa_{\ln y}d_{\theta'}$.
Indeed, the matrix $D_\theta$ changes by an overall sign
if $\theta$ changes by $2\pi$ and so does the matrix $G=B_{x}A_{\ln y}D_\theta$. 
Therefore we can find a unique $k\in\Z$ such that $\theta':= 2k\pi+\theta\in [0,2\pi)$
to have $d_{\theta'}=d_\theta$. This implies the first equality in (b). The latter follows from 
$d_{\theta+2k\pi}=d_\theta$ for all $\theta\in [0,2\pi)$ and $k\in\Z$.
{\hfill$\Box$}\bigskip

\noindent
{\bf Proof of Theorem \ref{fdPSL}.} First, 
denote
\[\hat F=\{G=B_x A_{\ln y} D_\theta \in\SL(2,\R) : x+iy\in F, \theta\in [0, 2\pi)\}.\] 
It is easy to see that $\hat F$ is open in $\SL(2,\R)$ and since the projection
$\pi: \SL(2,\R)\rightarrow \PSL(2,\R)$ is an open map and $\pi(\hat F)=\F$, it
follows that ${\cal F}$ is open in $\PSL(2,\R)$ (note that $d_{\theta+2k\pi}=d_\theta$ for $\theta\in [0,2\pi)$ and $k\in\Z$). 
To establish part (a) in the definition, we first claim that the closure of $\F$ in $\PSL(2,\R)$ is
\[\overline \F=
\{g=b_x a_{\ln y}d_\theta\,\in\PSL(2,\R) :\, x+iy\in \overline F,  \theta\in [0, 2\pi) \}.\]
Indeed, it suffices to check that
\begin{equation}\label{hatF} 
\overline {\hat F}=\{G=B_x A_{\ln y}D_\theta\, :\, x+iy\in \overline F, \theta\in [0, 2\pi)\},
\end{equation} 
where $\overline F$ denotes the closure of $F$ in $\H^2$
and $\overline{\hat F}$ denotes the closure of $\hat F$ in $\SL(2,\R)$.  
The set in the right-hand side of \eqref{hatF} is denoted by $cl(\hat F)$. 
For every $G\in cl\hat F$, we show that $G$ is the limit for some sequence $(G_n)\subset \hat F$. 
Writing $G=B_x A_{\ln y} D_\theta$, we have $x+iy\in \overline F$ by the definition of $cl(\hat F)$.  
Let $(x_n+iy_n)_n\subset  F$ be such that
$x_n+iy_n\rightarrow x+iy $ in $\H^2$ as $n\rightarrow \infty$. Then $x_n\rightarrow x$ 
as well as $y_n\rightarrow y$ in $\R$. Taking
$G_n = B_{x_n}A_{\ln y_n}D_\theta\in cl(\hat F)$, we obtain 
$G_n\rightarrow G$ in $\SL(2,\R)$ after a short check.

Next, for any $g=b_xa_{\ln y}d_\theta \in \PSL(2,\R)$ we have
$z:=x+iy\in \Phi(\gamma)(\overline F)$ for some $\gamma\in\Gamma$
as $F\subset \H^2$ is a fundamental domain for $\Gamma$.  
Take $\tilde z=\tilde x+i\tilde y=\Phi(\gamma^{-1})(z)\in\overline F$ and write $\gamma=[T]$ with 
$T=\scriptsize\Big(\begin{array}{cc}t_{11}&t_{12}\\t_{21}&t_{22}\end{array} \Big)
\in\SL(2,\R)$. Let 
\[\tilde \theta=\theta+2\arg(t_{21}\tilde z+t_{22})+2k\pi\in[0,2\pi)\]
for a unique $k\in\Z$ to obtain $h:=b_{\tilde x}a_{\ln \tilde y}d_{\tilde\theta}\in\overline{\F}$.
Thus 
\[ x+iy=\frac{t_{11}\tilde z +t_{12}}{t_{21}\tilde z+t_{22}} 
\quad \mbox{and}\quad \theta =\tilde \theta-2\arg(t_{21}\tilde z+t_{22})-2k\pi\]
imply $g=\gamma h\in\gamma\overline{\F}$ after a short
calculation. This completes the proof for (a) in the definition.

For part (b), suppose on the contrary that there exists 
$g\in \F \cap \gamma\F$ for some $\gamma\in\Gamma\setminus\{e\}$.
Then 
$g=b_xa_{\ln y}d_\theta$ and $g=\gamma b_{x'}a_{\ln y'} d_{\theta'}$
for $x+iy\in F$ and $x'+iy'\in F$. 
A short calculation shows that
$x+iy=\Phi(\gamma)(x'+iy')\in F \cap \Phi(\gamma)(F)$, which however contradicts the fact that $F$ is a fundamental domain. Thus $\F\cap \gamma \F=\varnothing$ for all $\gamma\in\Gamma\setminus\{e\}$.

Conversely, assume that $\F$ is a fundamental domain for $\Gamma$. Then $F\subset \H^2$ is open
since $\F$ is open.  
For any $z=x+iy\in\H^2$, then $g=b_xa_{\ln y}\in\PSL(2,\R)=\cup_{\gamma\in\Gamma}\gamma \overline{\F}$
implies that
$g=\gamma h$ for some $\gamma\in\Gamma$ and $h\in\overline{\F}=\{b_xa_{\ln y}d_\theta, x+iy\in\overline F, 0\leq \theta\leq 2\pi\}$.
Write $h=b_{\tilde x}a_{\ln\tilde y}d_{\tilde\theta}$. Then $\tilde z=\tilde x+i\tilde y\in \overline F$ and $z=\Phi(\gamma)(\tilde z)$ yield $z\in \Phi(\gamma)(\overline F)$. This proves (a) in Definition \ref{funddo-def}.   
Finally, assume that $z=x+iy\in F$ and $z=\Phi(\gamma)(z')$ for some $\gamma\in\Gamma\setminus\{e\} $ and $z'=x'+iy'\in F$.
Then take $g=b_xa_{\ln y}d_{\pi}$ and $h=b_{x'}a_{\ln y'}d_\theta$ with $\theta=2\arg(h_{21}z+h_{22})+2k\pi$
for a unique $k\in\Z$ such that $\theta \in [0,2\pi)$; here $h=\pi(H)$, $H=\scriptsize\Big(\begin{array}{cc}h_{11}&h_{12}\\h_{21}&h_{22}\end{array} \Big)$.  Then $g=\gamma h$ after a short computation.
This means that $\F\cap\gamma \F\ne \varnothing$, which is impossible since $\F$ is a fundamental domain.
{\hfill$\Box$}\bigskip

The next result follows directly from Proposition \ref{examfun} and Lemma \ref{fdPSL}. 
\begin{corollary}\label{fundaF0}
	(a) For $t>0$, the set
	\begin{equation}\label{et}
	\F_t=\{g=b_x a_{\ln y}d_\theta\in\PSL(2,\R) \,:\, x\in\R,
	1<y<e^t, \theta\in [0,2\pi)\}
	\end{equation}
	is a fundamental domain in $\PSL(2,\R)$ for the Fuchsian group $\langle a_t\rangle.$
	
	(b) For $t>0$, the set 
	\begin{equation}\label{et2}
	{\cal E}_t=\{g=b_x a_{\ln y}d_\theta\in\PSL(2,\R) \,:\,0< x<t,
	y>0, \theta\in [0,2\pi)\}
	\end{equation} is a fundamental domain in $\PSL(2,\R)$ for the Fuchsian group
	 $\langle b_t\rangle$.
	
\end{corollary}

It is well-known that $\PSL(2,\Z)$ is a Fuchsian group and the set
\[F=\Big\{z\in\H^2: |z|>1, |{\rm Re\,} z|<\frac{1}{2} \Big\} \]
is a fundamental domain of $\PSL(2,\Z)$ in $\H^2$ (see \cite[Proposition 9.18]{einsward}).
The following result follows from Lemma \ref{lemnak} and Theorem \ref{fdPSL}.
\begin{corollary}
	The set 
	\begin{eqnarray*}
	{\cal F}&=&\Big\{g=[G]\in\PSL(2,\R), G=\scriptsize\Big(\begin{array}{cc} a&b\\c&d \end{array}\Big)\in\SL(2,\R):\\
	&&	\ \ \ 2|ac+bd|<{c^2+d^2}, (ac+bd)^2+1> (c^2+d^2)^2 \Big\}
	\end{eqnarray*}	is a fundamental domain in $\PSL(2,\R)$ for $\PSL(2,\Z)$.
\end{corollary}

The next result shows us how to find a fundamental domain for a cyclic group 
as we know a fundamental domain for the cyclic group generated by
a conjugate element of its generator.

\begin{lemma}\label{fundaconj}
	Let $g_1$ and $g_2$ be conjugate in $\PSL(2,\R)$ and $g_2=h g_1 h^{-1}$ for $h\in\PSL(2,\R)$.
	Then if $\F_1\subset \PSL(2,\R)$ is a fundamental domain for
	$\langle g_1 \rangle $
	then $\F_2=h \F_1$ is a fundamental domain for
	$\langle g_2 \rangle$.  
\end{lemma}
\noindent
{\bf Proof.} Obviously $\overline{ \F_2}=h\overline{\F_1}$. 
Since $\F_1$ is a fundamental domain for $\langle g_1\rangle$, we have
\begin{eqnarray*}
	\bigcup_{j\in\Z} g_2^j \overline{\F_2}
	=\bigcup_{j\in\Z} hg_1^jh^{-1}h\overline{\F_1}
	=\bigcup_{j\in\Z} hg_1^j\overline{\F_1}
	= h(\bigcup_{j\in\Z} g_1^j\overline{\F_1})
	=h \PSL(2,\R)=\PSL(2,\R),
\end{eqnarray*}
and if $j\in\Z, g_2^j\ne e$, then $g_1^j\ne e$ yields  
\begin{eqnarray*}
	\F_2\cap g_2^j \F_2
	=h \F_1\cap hg_1^jh^{-1}h \F_1
	=h\F_1\cap hg_1^j \F_1
	=h(\F_1\cap g_1^j \F_1)=\varnothing.
\end{eqnarray*}
Also both $\F_1\subset \PSL(2,\R)$
and $\F_2\subset \PSL(2,\R)$ are open. 

{\hfill$\square$}\bigskip

Recall that every hyperbolic (resp. parabolic) element is conjugate with  $a_t$ (resp. $b_t$) for some $t\in\R$.  
Note that $\langle a_t\rangle=\langle a_{-t}\rangle$ and $\langle b_t\rangle=\langle b_{-t}\rangle$.
The next result follows from the preceding lemma.
\begin{proposition} Let $g\in \PSL(2,\R)$ be a hyperbolic element
	(resp. parabolic element). If  $h\in\PSL(2,\R)$ and $t\in\R$ are such that $g=h^{-1}a_th$
	(resp. $g=h^{-1}b_t h$) then ${\cal F}=h{\cal F}_{|t|}$
	(resp. ${\cal E}=h{\cal E}_{|t|}$) is a fundamental domain for 
	$\Gamma=\langle g\rangle $,  where  ${\cal F}_{|t|}$ (resp. ${\cal E}_{|t|}$) is a fundamental domain 
	for $\langle a_t\rangle$  (resp. $\langle b_t\rangle$) given by \eqref{et} (resp. \eqref{et2}).
\end{proposition}

Next we define $\Theta:T^1\H^2\to \PSL(2,\R)$ by $\Theta(z,\xi)=g$ for $(z,\xi)\in T^1\H^2$,
where $g\in\PSL(2,\R)$ satisfies $\D g(i,i)=(z,\xi)$. 
Then $\Theta$ is well-defined and bijective owing to Lemma \ref{unig}.
Note that there exist metrics on $\PSL(2,\R)$ and $T^1\H^2$ such that $\Theta$ is an isometry.
\begin{lemma}\label{cc}
	Let $F\subset\H^2$ and denote $T^1F=\{(z,\xi)\in T^1\H^2: z\in F \}$.
	Then  \[\Theta(T^1F)=\{g\in\PSL(2,\R):g=b_x a_{\ln y}d_\theta: x+iy\in F,\theta\in [0,2\pi) \}.\]
\end{lemma} 
\noindent
{\bf Proof\,:} For any $g=b_x a_{\ln y}d_\theta\in\PSL(2,\R)$
with $x+iy\in F$, if $g=\big\{\pm \big({\scriptsize\begin{array}{cc}a&b\\c&d \end{array}}\big) \big\}$ then
we take $z=x+iy\in F$ and $\xi =\Re \xi +i\Im \xi$ satisfying 
\[ x=\frac{ac+bd}{c^2+d^2},\ y=\frac{1}{c^2+d^2},\ {\rm Re\,}\xi=\frac{2cd}{(c^2+d^2)^2},
\ \Im\xi=\frac{d^2-c^2}{(c^2+d^2)^2}. \]
Then $\|\xi\|_z=\frac{|\xi|}{y}=1$ means that $(z,\xi)\in T^1F$
and  \eqref{dg} shows $\Theta(z,\xi)=g.$
On the other hand, for $(z,\xi)\in T^1F$ and $\Theta(z,\xi)=g\in\PSL(2,\R)$. If $g=\Big\{\pm \scriptsize\Big(\begin{array}{cc}a&b\\c&d \end{array}\Big)\Big\}=b_xa_{\ln y}d_\theta$
then  $ x=\frac{ac+bd}{c^2+d^2}, y=\frac{1}{c^2+d^2}$ by Lemma \ref{lemnak} (a).
Once again \eqref{dg} implies that $z=x+iy\in F$. This completes the proof.
{\hfill $\Box$} 

The relation of fundamental domains in $\H^2$ and in $T^1\H^2$ is the following:

\begin{theorem} Let $\Gamma$ be a Fuchsian group. 
A set	 $F\subset\H^2$ is a fundamental domain for $\Gamma$ if and only if $T^1F\subset T^1\H^2$ is a 
	fundamental domain for $\Gamma$. 
\end{theorem}
\noindent{\bf Proof\,:}
Let $F\subset \H^2$ and ${\cal F}=\{ g=b_xa_{\ln y}d_\theta, x+iy\in F, \theta\in [0,2\pi)\}\subset \PSL(2,\R)$.
Then $\Theta^{-1}({\cal F})= T^1F$ by Lemma \ref{cc} and this follows from Theorem \ref{fdPSL}
and the fact that $\Theta$ is an isometry.
{\hfill$\Box$}

\bibliographystyle{mystyle}

\noindent	{\bf Acknowledgments:} This work is supported by Vietnam National Foundation
for Science and Technology Development (Grant No. 101.02-2020.21).

\end{document}